\documentclass[10pt,a4paper,reqno]{amsart}
\usepackage{amsthm}
\usepackage{amsmath}
\usepackage{amsmath, amssymb}
\usepackage{type1cm}
\usepackage{cite}
\usepackage{bm}
\usepackage{cases}
\usepackage{here}
\usepackage[dvipdfmx]{graphicx}
\theoremstyle{definition}
\newtheorem{defi}{Definition}
\newtheorem{prop}{Proposition}

\newtheorem{theorem}{Theorem}
\newtheorem{corollary}{Corollary}
\newcommand{\eqdef}{\stackrel{\mathrm{def}}{=}}
\DeclareMathOperator*{\esssup}{ess\,sup}
\begin{document} 
\title{$L^p$-norm inequality using q-moment and its applications}
\author{Tomohiro Nishiyama}
\begin{abstract}
For a measurable function on a set which has a finite measure, an inequality holds between two Lp-norms. In this paper, we show similar inequalities for the Euclidean space and the Lebesgue measure by using a q-moment which is a moment of an escort distribution. As applications of these inequalities, we first derive upper bounds for the Renyi and the Tsallis entropies with given q-moment and derive an inequality between two Renyi entropies. Second, we derive an upper bound for the probability of a subset in the Euclidean space with given Lp-norm on the same set. 
\\

\end{abstract}
\date{}
\maketitle
\bibliographystyle{plain}
\section{Introduction}
We consider the measure space $(\mathbb{R}^n, \mathcal{M}, \mu)$, where $\mathcal{M}$ is a $\sigma$-algebra and $\mu$ is the Lebesgue measure.
When $S\in\mathcal{M}$ has a finite measure, for $0<q<r$ and measurable function $f$ on $S$, the following inequality holds.
\begin{align}
\label{eq_lp_inequality}
\|f\|_q\leq \mu(S)^{\frac{1}{q}-\frac{1}{r}}\|f\|_r,
\end{align}
where $\|f\|_p$ is a $L^p$-norm\cite{rudin2006real} defined as follows.

\begin{defi}
Let $S\in\mathcal{M}$.
For $0<p<\infty$, 
\begin{align}
\|f\|_p\eqdef {\bigl(\int_{S} |f(x)|^p \mathrm{d}^nx\bigr)}^{\frac{1}{p}}. \nonumber
\end{align}
For $p=\infty$,
\begin{align}
\|f\|_\infty\eqdef \esssup_{x\in S} |f|. \nonumber
\end{align}
\end{defi}

In this paper, we show a similar inequality between two $L^p$-norms by using a $q$-moment when $S=\mathbb{R}^n$ (See Theorem 1 and 2 in section 2).

For $0<q<r\leq\infty$ and $n=1$,
\begin{align}
\label{eq_new_inequality}
\|f\|_q\leq {(C{\mu_{q,\alpha}}^{\frac{1}{\alpha}})}^{\frac{1}{q}-\frac{1}{r}}\|f\|_r,
\end{align} 
where $C$ is a constant and $\mu_{q,\alpha}$ is the $\alpha$-th order $q$-moment defined as follows. .

\begin{defi}
Let $f:\mathbb{R}^n\rightarrow\mathbb{R}$ be a measurable function which satisfies $0<\|f\|_q < \infty$.
Let $X\in \mathbb{R}^n$.
We define a q-expectation value \cite{tsallis1998role,abe2005necessity}.
\begin{align}
E_q[X]\eqdef \frac{\int_{\mathbb{R}^n} x|f(x)|^q \mathrm{d}^n x}{\int_{\mathbb{R}^n} |f(x)|^q \mathrm{d}^n x} \nonumber
\end{align}
Especially, when $q=1$, we write $E_1[X]$ as $E[X]$.
\end{defi}
 $E_q[X]$ is also a expected value of a escort distribution $\frac{|f(x)|^q}{\int_{\mathbb{R}^n} |f(x)|^q \mathrm{d}^n x}$\cite{beck1995thermodynamics}.

\begin{defi}
\label{def_q_moment}
Let $f:\mathbb{R}\rightarrow\mathbb{R}$ be a measurable function and $X\in \mathbb{R}$.
For $b\in \mathbb{R}$ and $\alpha >0$, we define the $\alpha$-th order $q$-moment as follows.
\begin{align}
\mu_{q,\alpha}\eqdef E_q[|X-b|^\alpha] \nonumber
\end{align}
When $b=E_q[X]$, $\mu_{q,\alpha}$ is the central $q$-moment.
\end{defi}
$\mu_{q,\alpha}$ is also the $\alpha$-th order moment of a escort distribution.

In (\ref{eq_new_inequality}), ${\mu_{q,\alpha}}^{\frac{1}{\alpha}}$ corresponds to $\mu(S)$ in (\ref{eq_lp_inequality}) and we can interpret ${\mu_{q,\alpha}}^{\frac{1}{\alpha}}$ as the ``range'' of the region function $f$ spreads.

As applications of the inequality (\ref{eq_new_inequality}) and the multivariate version of  (\ref{eq_new_inequality}), we derive an inequality between two R\'{e}nyi entropies\cite{renyi1961measures} and derive upper bounds for the R\'{e}nyi and the Tsallis entropies with given $q$-moment. We also obtain an upper bound for the Shannon differential entropy \cite{conrad2004probability} as a limit of the R\'{e}nyi entropy.

Furthermore, we derive an upper bound for the probability of a subset in $\mathbb{R}^n$ with given $L^p$-norm on the same set. This is a generalization of the result in\cite{nishiyama2018improved}.

\section{Main Results}
\begin{theorem}
\label{th_1dim_norm_inequality}
Let $f:\mathbb{R}\rightarrow \mathbb{R}$ be a measurable function with finite $q$-moment. Let $\|f\|_q,\|f\|_r<\infty$.

For $0<q<r\leq\infty$,
\begin{align}
\|f\|_q\leq {(C{\mu_{q,\alpha}}^{\frac{1}{\alpha}})}^{\frac{1}{q}-\frac{1}{r}}\|f\|_r, 
\end{align} 
where $C$ is a constant which only depends on $\alpha$.
The example value of $C$ is $C=\frac{2}{\alpha}\Gamma(\frac{1}{\alpha})(\alpha e)^{\frac{1}{\alpha}}$.
\end{theorem}
\noindent\textbf{Proof.}
For a non-negative convex function $\phi_t$, we consider the following value.
\begin{align}
\label{eq_basic}
V=\int_{\mathbb{R}} |f(x)|^q\phi_t(\mu_{q,\alpha}^{-1}|x-b|^\alpha) \mathrm{d}x, 
\end{align}
where the function $\phi_t(x)$ satisfies $\phi_t(1)=1$.
We transform this equation as follows.
\begin{align}
V=\int_{\mathbb{R}} |f(x)|^q \mathrm{d}x \times E_q[\phi_t(\mu_{q,\alpha}^{-1}|X-b|^\alpha)]
=\|f\|_q^q E_q[\phi_t(\mu_{q,\alpha}^{-1}|X-b|^\alpha)]
\end{align}
Applying the Jensen's inequality to this equation and using Definition \ref{def_q_moment} give
\begin{align}
\label{eq_jensen}
V\geq \|f\|_q^q\phi_t(\mu_{q,\alpha}^{-1}E_q[|X-b|^\alpha]) =\|f\|_q^q\phi_t(1)=\|f\|_q^q.
\end{align}

Furthermore, for $1\leq s, t \leq \infty$, by applying the H\"{o}lder's inequality to (\ref{eq_basic}), we have
\begin{align}
V\leq \||f|^q\|_s{(\int_{\mathbb{R}} \phi_t(\mu_{q,\alpha}^{-1}|x-b|^\alpha)^t \mathrm{d}x)}^{\frac{1}{t}},
\end{align}
where $\frac{1}{s}+\frac{1}{t}=1$.
By assumption $q<r$, we can put $r=qs$ and $1\leq t<\infty$. 
Then, we have
\begin{align}
\label{eq_holder}
V\leq \|f\|_r^q{(\int_{\mathbb{R}} \phi_t(\mu_{q,\alpha}^{-1}|x-b|^\alpha)^t \mathrm{d}x)}^{\frac{1}{t}},
\end{align}
where we use $\||f|^q\|_s=\|f\|_{qs}^q=\|f\|_r^q$.

Since $\phi_t(x)=\exp\bigl(-\frac{\beta}{t}(x-1)\bigr)$ is a convex function and satisfies $\phi_t(1)=1$, substituting $\phi_t(x)=\exp\bigl(-\frac{\beta}{t}(x-1)\bigr)$ into RHS of this inequality, we have
\begin{align}
\int_{\mathbb{R}} \phi_t(\mu_{q,\alpha}^{-1}|x-b|^\alpha)^t \mathrm{d}x=\int_{\mathbb{R}}\exp\bigl(-\beta(\mu_{q,\alpha}^{-1}|x-b|^\alpha-1)\bigr)\mathrm{d}x,
\end{align}
where $\beta > 0$.
Changing from the variable $x$ to $y=\mu_{q,\alpha}^{-\frac{1}{\alpha}}(x-b)$ gives
\begin{align}
\label{eq_rhs}
\int_{\mathbb{R}} \phi_t(\mu_{q,\alpha}^{-1}|x-b|^\alpha)^t \mathrm{d}x&=2\exp(\beta){\mu_{q,\alpha}}^{\frac{1}{\alpha}
}\int_0^\infty\exp\bigl(-\beta y^\alpha)\mathrm{d}y
&=\frac{2}{\alpha}\Gamma(\frac{1}{\alpha}){\mu_{q,\alpha}}^{\frac{1}{\alpha}}\exp(\beta)\beta^{-\frac{1}{\alpha}}.
\end{align}
$\exp(\beta)\beta^{-\frac{1}{\alpha}}$ has a minimum value at $\beta=\frac{1}{\alpha}$.
Substituting this condition into (\ref{eq_rhs}) and (\ref{eq_holder}) gives
\begin{align}
\label{eq_holder_2}
V\leq \|f\|_r^q{(C{\mu_{q,\alpha}}^{\frac{1}{\alpha}})}^{\frac{1}{t}},
\end{align}
where $C=\frac{2}{\alpha}\Gamma(\frac{1}{\alpha})(\alpha e)^{\frac{1}{\alpha}}$.
Combining (\ref{eq_jensen}) and (\ref{eq_holder_2}) gives
\begin{align}
\label{eq_final}
\|f\|_q\leq \|f\|_r{(C{\mu_{q,\alpha}}^{\frac{1}{\alpha}})}^{\frac{1}{qt}}.
\end{align}
Combining  $\frac{1}{s}+\frac{1}{t}=1$ and $r=qs$, we have $\frac{1}{qt}=\frac{1}{q}-\frac{1}{r}$.
Substituting this equation into (\ref{eq_final}), the result follows.

Next, we prove the multivariate version.
\begin{defi}
Let $f:\mathbb{R}^n\rightarrow\mathbb{R}$ be a measurable function.

For $b\in \mathbb{R}^n$ and $X\in \mathbb{R}^n$, we define a multivariate $q$-moment as follows.
\begin{align}
\Sigma_{q,b}\eqdef E_q[(X-b)(X-b)^T], \nonumber
\end{align}
where $T$ denotes the transpose of a vector.

When $b=E_q[X]$, we write $\Sigma_{q,b}$ as $\Sigma_q$ and $\Sigma_q$ is equal to a $q$-covariance matrix.

When $q=1$, we write $\Sigma_{q,b}$ as $\Sigma_b$ and $\Sigma$ denotes a covariance matrix.
\end{defi}
\begin{theorem}
\label{th_multi_norm_inequality}
Let $f:\mathbb{R}^n\rightarrow\mathbb{R}$ be a measurable function with finite multivariate $q$-moment. Let $\|f\|_q,\|f\|_r<\infty$.

For $0<q<r\leq\infty$ and $n\geq 1$
\begin{align}
\|f\|_q\leq {(C{(\det\Sigma_{q,b})}^{\frac{1}{2}})}^{\frac{1}{q}-\frac{1}{r}}\|f\|_r. 
\end{align} 
$C$ is a constant which only depends on $n$.
The example value of $C$ is $C={(2\pi e)}^{\frac{n}{2}}$.
\end{theorem}
\noindent\textbf{Proof.}
We can prove this theorem in the same way as the theorem \ref{th_1dim_norm_inequality}.

First, we consider the following value.
\begin{align}
\label{eq_multi_basic}
V=\int_{\mathbb{R}^n} |f(x)|^q\phi_t((x-b)^T\Sigma_{q,b}^{-1}(x-b)) \mathrm{d}^nx, 
\end{align}
where $\phi_t$ is a non-negative convex function which satisfies $\phi_t(n)=1$.
By applying the Jensen's inequality to this equation in the same way as Theorem \ref{th_1dim_norm_inequality}, we get
\begin{align}
\label{eq_multi_jensen}
V\geq \|f\|_q^q\phi_t(n)=\|f\|_q^q.
\end{align}
Next, for $1\leq s, t \leq \infty$, by applying the H\"{o}lder's inequality to (\ref{eq_multi_basic}) and putting $r=qs$, we have
\begin{align}
\label{eq_multi_holder_1}
V\leq \|f\|_r^q{\biggl(\int_{\mathbb{R}^n} \phi_t\bigl((x-b)^T\Sigma_{q,b}^{-1}(x-b)\bigr)^t \mathrm{d}x\biggr)}^{\frac{1}{t}}.
\end{align}
Substituting $\phi_t(x)=\exp\bigl(-\frac{\beta}{t}(x-n)\bigr)$ into this inequality gives
\begin{align}
\label{eq_multi_integral_1}
\int_{\mathbb{R}^n} \phi_t\bigl((x-b)^T\Sigma_{q,b}^{-1}(x-b)\bigr)^t \mathrm{d}x=\exp(\beta n)\int_{\mathbb{R}^n}\exp\bigl(-\beta(x-b)^T\Sigma_{q,b}^{-1}(x-b)\bigr)\mathrm{d}^nx.
\end{align}
Changing the variable from $x$ to $y=\Sigma_{q,b}^{-\frac{1}{2}}(x-b)$ gives
\begin{align}
\label{eq_multi_rhs}
\exp(\beta n)\int_{\mathbb{R}^n}\exp\bigl(-\beta(x-b)^T\Sigma_{q,b}^{-1}(x-b)\bigr)\mathrm{d}^nx=\exp(\beta n){(\det\Sigma_{q,b})}^{\frac{1}{2}}{\biggl(\frac{\pi}{\beta}\biggr)}^{\frac{n}{2}}.
\end{align}
$\exp(\beta)\beta^{-\frac{1}{2}}$ has a minimum value at $\beta=\frac{1}{2}$.
Substituting this condition into (\ref{eq_multi_rhs}) and combining (\ref{eq_multi_holder_1}) and (\ref{eq_multi_integral_1}) give
\begin{align}
\label{eq_multi_holder_2}
V\leq \|f\|_r^q{(C{(\det\Sigma_{q,b})}^{\frac{1}{2}})}^{\frac{1}{t}},
\end{align}
where $C={(2\pi e )}^{\frac{n}{2}}$.
Combining (\ref{eq_multi_jensen}) and (\ref{eq_multi_holder_2}) gives
\begin{align}
\label{eq_multi_final}
\|f\|_q\leq \|f\|_r{(C{(\det\Sigma_{q,b})}^{\frac{1}{2}})}^{\frac{1}{qt}}.
\end{align}
Combining  $\frac{1}{s}+\frac{1}{t}=1$ and $r=qs$, we have $\frac{1}{qt}=\frac{1}{q}-\frac{1}{r}$.
Substituting this equation into (\ref{eq_multi_final}), the result follows.

\section{Applications}
In this section, $C$ denotes a constant which only depends on the dimension $n$.
The example value of $C$ is $(2\pi e)^{\frac{n}{2}}$.
\subsection{Application for the R\'{e}nyi and the Tsallis entropies}
We derive upper bounds for the R\'{e}nyi and the Tsallis entropies with given $q$-covariance matrix and derive an inequality between two R\'{e}nyi entropies by using Theorem \ref{th_multi_norm_inequality}.

\begin{corollary}
\label{co_renyi_variance_inequality}
Let $f$ be a probability density function on $\mathbb{R}^n$ with finite $q$-covariance matrix and $h_p(X)\eqdef \frac{p}{1-p}\log \|f\|_p$ be the R\'{e}nyi entropy.

For $p>1$,
\begin{align}
h_p(X)\leq  \frac{1}{2}\log({\det\Sigma})+\log C.
\end{align}

For $0<p<1$, 
\begin{align}
h_p(X)\leq  \frac{1}{2}\log({\det\Sigma_p})+\log C.
\end{align}
\end{corollary}  
\noindent\textbf{Proof.}
For $p>1$, by putting $q=1$, $r=p$, $b=E_q[X]$ and using $\|f\|_1=1$ in Theorem \ref{th_multi_norm_inequality}, the result follows.

For $0<p<1$, by putting $r=1$, $q=p$, $b=E_q[X]$ and using $\|f\|_1=1$ in Theorem \ref{th_multi_norm_inequality}, the result follows.

We can derive the optimal constant $C$ by using the distribution that maximizes the R\'{e}nyi entropy\cite{johnson2007some}.

In the limit $p\rightarrow 1$, the R\'{e}nyi entropy is the Shannon differential entropy and we have 
\begin{align}
h(X)\leq  \frac{1}{2}\log({\det\Sigma})+\log C,
\end{align}
where $h(X)$ is the Shannon differential entropy.
When $C=(2\pi e)^{\frac{n}{2}}$, this inequality is consistent with the well-known upper bound of the Shannon entropy.

\begin{corollary}
Let $f$ be a probability density function on $\mathbb{R}^n$ with finite $q$-covariance matrix and $h_p(X)\eqdef \frac{p}{1-p}\log \|f\|_p$ be the R\'{e}nyi entropy.

For $0<q<r\leq \infty$,
\begin{align}
\frac{1-q}{q}h_q(X)\leq  \frac{1-r}{r}h_r(X)+\biggl(\frac{1}{q}-\frac{1}{r}\biggr)\bigl(\frac{1}{2}\log(\det\Sigma_q)+\log C\bigr).
\end{align}

\end{corollary}  
\noindent\textbf{Proof.}
Taking the logarithm of Theorem \ref{th_multi_norm_inequality}, the result follows.

\begin{corollary}
Let $f$ be a probability density function on $\mathbb{R}^n$ with finite $q$-covariance matrix and $S_q(X)\eqdef \frac{1}{q-1}(1-\|f\|_q^q)$ be the Tsallis entropy\cite{tsallis1988possible}.

For $q>1$,
\begin{align}
\exp_q(S_q(X))\leq  C{(\det\Sigma)}^{\frac{1}{2}},
\end{align}
where $\exp_q(x)\eqdef [1+(1-q)x]^{\frac{1}{1-q}}$.

For $0<q<1$, 
\begin{align}
\exp_q(S_q(X))\leq  C{(\det\Sigma_q)}^{\frac{1}{2}}.
\end{align}
\end{corollary}  
\noindent\textbf{Proof.}
From the definition of the Tsallis entropy and $\exp_q(x)$, we have $\exp_q(S_q(X))=\|f\|_q^{\frac{q}{1-q}}$.
The rest part of the proof is almost the same as Corollary \ref{co_renyi_variance_inequality}.

\subsection{Application for the upper bound of a probability}
We derive an upper bound for the probability of a subset with given $L^p$-norm on the same set. 
Since we use some functions in this subsection, we use the following notation. \\
\noindent\textbf{Notation.}
\begin{itemize}
\item
For a non-negative measurable function $F$,
\begin{align}
E_F[X]\eqdef \frac{\int_{\mathbb{R}^n}  x F(x)\mathrm{d}^nx}{\int_{\mathbb{R}^n} F(x) \mathrm{d}^nx}. \nonumber \\ 
\Sigma_{F,b}\eqdef E_F[(X-b)(X-b)^T].\nonumber
\end{align}
When $b=E_q[X]$, we write $\Sigma_{F,b}$ as $\Sigma_F$.

\item
For $\Omega\subseteq\mathbb{R}^n$, 
\[
  I_\Omega(x) = \begin{cases}
    1 & (x\in \Omega) \\
    0 & (x\notin \Omega)
  \end{cases}
\]
\end{itemize}
\begin{prop}
Let $f$ be a probability density function on $\mathbb{R}^n$ with finite covariance matrix $\Sigma_f$.
Let $\Omega\subseteq\mathbb{R}^n$ and $P(\Omega)$ be a probability of $\Omega$.

For $r>1$ and $n\geq 1$,
\begin{align}
\label{eq_prob_upper_bound}
P(\Omega)^{1+\frac{n}{2}-\frac{n}{2r}}\leq  {(C{(\det\Sigma_f)}^{\frac{1}{2}})}^{1-\frac{1}{r}}\|fI_\Omega\|_r.
\end{align}
\end{prop}
\noindent\textbf{Proof.}
First, we put 
\begin{align}
 \psi(x)\eqdef {(\Sigma_f)}^{-\frac{1}{2}}(x-E_f[X]) \\ \nonumber
y=\psi(x) \\ \nonumber
\hat{\Omega}=\psi(\Omega)\\ \nonumber
\hat{f}(y)\eqdef{(\det\Sigma_f)}^{\frac{1}{2}}f(x) \\ \nonumber
g(x)\eqdef f(x)I_\Omega(x) \\ \nonumber
\hat{g}(y)\eqdef \hat{f}(y)I_{\hat{\Omega}}(y) .
\end{align}
Since ${(\det\Sigma_f)}^{\frac{1}{2}}\mathrm{d}^ny=\mathrm{d}^nx$, we have
\begin{align}
\int_{\mathbb{R}^n} \hat{g}(y)^r \mathrm{d}^ny&=\int_{\hat{\Omega}} \hat{f}(y)^r \mathrm{d}^ny \\ \nonumber
&={(\det\Sigma_f)}^{\frac{r-1}{2}}\int_{\Omega} f(x)^r\mathrm{d}^nx \\ \nonumber
&={(\det\Sigma_f)}^{\frac{r-1}{2}}\int_{\mathbb{R}^n} g(x)^r\mathrm{d}^nx.
\end{align}

From this equation, we have
\begin{align}
\label{eq_change_variable}
\|\hat{g}\|_r={(\det\Sigma_f)}^{\frac{r-1}{2r}}\|g\|_r={(\det\Sigma_f)}^{\frac{1}{2}(1-\frac{1}{r})}\|g\|_r.
\end{align}
Especially, when $r=1$, $\|\hat{g}\|_1=\|g\|_1=P(\Omega)$.

Furthermore, using the GM-AM inequality ${(\det A)}^{\frac{1}{n}} \leq \frac{1}{n}\mathrm{Tr}A $, we have
\begin{align}
\label{eq_trace_inequality}
\det\Sigma_{\hat{g},0}\leq  {\bigl(\frac{1}{n}\mathrm{Tr}\Sigma_{\hat{g},0}\bigr)}^n
= \frac{1}{\|\hat{g}\|_1^n}{\bigl(\frac{1}{n}\sum_i \int_{\mathbb{R}^n} y_i^2 \hat{g}(y) \mathrm{d}^ny\bigr)}^n \\ \nonumber
= \frac{1}{\|g\|_1^n}{\bigl(\frac{1}{n}\sum_i \int_{\mathbb{R}^n} y_i^2 \hat{g}(y) \mathrm{d}^ny\bigr)}^n,
\end{align}
and
\begin{align}
\label{eq_trace_inequality2}
\sum_i \int_{\mathbb{R}^n} y_i^2 \hat{g}(y) \mathrm{d}^ny\leq \sum_i \int_{\mathbb{R}^n} y_i^2 \hat{f}(y) \mathrm{d}^ny\\ \nonumber
=E_{\hat{f}}[\mathrm{Tr}(YY^T)] \\ \nonumber
=E_{f}[\mathrm{Tr}(\Sigma_f^{-1}(X-E_f[X])(X-E_f[X])^T)] \\ \nonumber
=\mathrm{Tr}(\Sigma_f^{-1}E_{f}[(X-E_f[X])(X-E_f[X])^T])=n. 
\end{align}
(\ref{eq_trace_inequality}) and (\ref{eq_trace_inequality2}) give
\begin{align}
\label{eq_trace_inequality3}
\det\Sigma_{\hat{g},0}\leq  \frac{1}{\|g\|_1^n}.
\end{align}

Applying Theorem \ref{th_multi_norm_inequality} for $\hat{g}$, $q=1$ and $b=0$ gives
\begin{align}
\|\hat{g}\|_1\leq {(C{(\det\Sigma_{\hat{g},0})}^{\frac{1}{2}})}^{1-\frac{1}{r}}\|\hat{g}\|_r.
\end{align}
By using (\ref{eq_change_variable}) and (\ref{eq_trace_inequality3}), we have
\begin{align}
\|g\|_1\leq {(C{(\det\Sigma_{\hat{g},0})}^{\frac{1}{2}})}^{1-\frac{1}{r}}\|\hat{g}\|_r\leq \frac{1}{\|g\|_1^{\frac{n}{2}(1-\frac{1}{r})}}{(C{(\det\Sigma_f)}^{\frac{1}{2}})}^{1-\frac{1}{r}}\|g\|_r.
\end{align}
By transforming this equation and using $P(\Omega)=\|g\|_1$, the result follows.

\begin{corollary}
\label{cor_prob_upperbound}
Let $f$ be a probability density function on $\mathbb{R}^n$ with finite covariance matrix $\Sigma_f$.
Let $\Omega\subseteq\mathbb{R}^n$ and $P(\Omega)$ be a probability of $\Omega$.

For $n\geq 1$,
\begin{align}
P(\Omega)\leq {\bigl(C{(\det\Sigma_f)}^{\frac{1}{2}}\|fI_\Omega\|_\infty\bigr)}^{\frac{2}{n+2}}.
\end{align}
\end{corollary}
\noindent\textbf{Proof.}
By substituting $r=\infty$ into (\ref{eq_prob_upper_bound}), the result follows.

Since $\|fI_\Omega \|_\infty=\esssup_{x\in \Omega} f(x)$, when the supremum of $f(x)$ in $\Omega$ is given, we can derive the probability upper bound by using Corollary \ref{cor_prob_upperbound}.

 
\section{Conclusion}
In the first half, we have shown inequalities between two $L^p$-norms by using the $q$-moment for the Euclidean space and the Lebesgue measure.

In the latter half, by applying these inequalities to probability theory, we have derived the inequality that holds between two R\'{e}nyi entropies, and have derived upper bounds for the R\'{e}nyi and the Tsallis entropies with given $q$-moment. 
In particular, by using the result of the R\'{e}nyi entropy, we have shown an upper bound for the Shannon entropy in the limit $q\rightarrow 1$.

Furthermore, we have derived the upper bound for the probability of the subset in $\mathbb{R}^n$ with given $L^p$-norm on the same set. 

We hope we will find the optimal constants $C$ for each inequality.

\bibliography{reference_v1}

\begin{thebibliography}{1}

\bibitem{abe2005necessity}
Sumiyoshi Abe and GB~Bagci.
\newblock Necessity of q-expectation value in nonextensive statistical
  mechanics.
\newblock {\em Physical Review E}, 71(1):016139, 2005.

\bibitem{beck1995thermodynamics}
Christian Beck and Friedrich Sch{\"o}gl.
\newblock {\em Thermodynamics of chaotic systems: an introduction}.
\newblock Number~4. Cambridge University Press, 1995.

\bibitem{conrad2004probability}
Keith Conrad.
\newblock Probability distributions and maximum entropy.
\newblock {\em Entropy}, 6(452):10, 2004.

\bibitem{johnson2007some}
Oliver Johnson and Christophe Vignat.
\newblock Some results concerning maximum r{\'e}nyi entropy distributions.
\newblock In {\em Annales de l'Institut Henri Poincar{\'e} (B) Probability and
  Statistics}, volume~43, pages 339--351. No longer published by Elsevier,
  2007.

\bibitem{nishiyama2018improved}
Tomohiro Nishiyama.
\newblock Improved chebyshev inequality: new probability bounds with known
  supremum of pdf.
\newblock {\em arXiv preprint arXiv:1808.10770}, 2018.

\bibitem{renyi1961measures}
Alfr{\'e}d R{\'e}nyi.
\newblock On measures of entropy and information.
\newblock Technical report, HUNGARIAN ACADEMY OF SCIENCES Budapest Hungary,
  1961.

\bibitem{rudin2006real}
Walter Rudin.
\newblock {\em Real and complex analysis}.
\newblock Tata McGraw-Hill Education, 2006.

\bibitem{tsallis1988possible}
Constantino Tsallis.
\newblock Possible generalization of boltzmann-gibbs statistics.
\newblock {\em Journal of statistical physics}, 52(1-2):479--487, 1988.

\bibitem{tsallis1998role}
Constantino Tsallis, RenioS Mendes, and Anel~R Plastino.
\newblock The role of constraints within generalized nonextensive statistics.
\newblock {\em Physica A: Statistical Mechanics and its Applications},
  261(3-4):534--554, 1998.

\end{thebibliography}
\end{document}